\documentclass[leqno,final,12pt,titlepage]{article}
\usepackage{amsfonts}

\usepackage{proc2e}

\input{tcilatex}

\begin{document}

\title{Some Results on Algebraic and Geometric Characterization of Linear Systems
Models for Time Series Analysis }
\author{Jo\~{a}o Jos\'{e} de Farias Neto \\
Institute for Advanced Studies - IEAv-CTA \\
Divis\~{a}o de \ Inform\'{a}tica, joaojfn@ieav.cta.br \\
Rod. dos Tamoios, km 5,5, Cep - 12.228-840 \\
S\~{a}o Jos\'{e} dos Campos, SP, Brazil}
\maketitle

\begin{abstract}
\begin{description}
\item  It is shown that in the multivariate case the orders p, of the AR
part, and q, of the MA part, are not invariants of the time series. Thus, it
is concluded that it only makes sense to define the class of ARMA(p,p)-
irreducible models, where p is the biggest of the system's Kronecker
indices. This class is shown not to be a differentiable manifold, but to
contain one, which is a generic subset of systems with all Kronecker indices
equal to p. A formula which gives the metric tensor for riemannian manifolds
of linear systems as a line integral in the complex plane is introduced for
deterministic and stochastic cases and some tensors are obtained with it.

\item  MSC-class: 93C05; 93A30; 93B29; 93B30
\end{description}
\end{abstract}

\section{Introduction}

Since the publication of Box and Jenkins book \cite{BJ}, the statistical \
methodology proposed by these authors has been disseminated and widely used
to build mathematical models of time series. These models are essentially
discrete time differential equations (thus, indeed difference equations), of
which the order and the parameters (supposed constant) are determined as a
function only of a unique sample of its output sequence - the time series -
which is considered as a particular realization of a gaussian stochastic
process resulting from a forced term (the input of the system) which is also
considered as an unobservable gaussian process. The orders p of the
difference equation and q of the input noise are identified with the aid of
correlation and partial correlation functions.

In the beginning of the 80's, Tiao and Tsay \cite{TT} , of the statistics
department of Wisconsin University, proposed a methodology analogous to Box
and Jenkin's to the multivariate case, that is, the one in which several
time series are treated simultaneously, being viewed as a vector time
series. The advantage of this approach is that the separate treatment of the
series would not take in consideration their interrelations, which otherwise
would allow the building of a model of greater predictive capacity. The
proposed methodology is based upon a direct generalization of the
correlation and partial correlation functions, to identify the orders p and
q of the model.

It happens, nevertheless, that there is not a method to determine isolately
the orders p and q of a multivariate ARMA(p,q)\ model nor is it possible to
come up with one, simply because, contrarily to the univariate case, a
multivariate ARMA model has not independent intrinsic p and q orders. This
will be shown using the polynomial representation of linear systems and
unimodular matrices.

The necessary correction is also proposed, with \ the introduction of the
ARMA(p,p)-irreducible models. Given the formal resemblance between this new
class of models and the class of linear systems of a common McMillan degree,
a natural question is if this new class also constitutes a differentiable
manifold. If this was the case, a whole new set of parametrizations would be
unveiled, as it happened with the so-called overlapping parametrizations.

Since the space of all linear dynamical systems with a common McMillan
degree and the generic set in the ARMA(p,p)-irreducible class are
differentiable manifolds, it makes sense to be able to do \textit{\
identification on a manifold }(see references \cite{HZ} and \cite{PT} for
the basic results on this area). The resulting path in the systems space
until convergence will be almost independent of the parametrization \cite{CH}%
. This could be termed \textit{coordinate-free identification.}

To use this approach, the Riemannian metric tensor G must be computed, since
the Riemannian gradient is given by G$^{-1}\nabla $J, that is, G$^{-1}$
times the gradient of a convenient objective function J. Furthermore, the
local geometric properties of the manifold are defined by G and, as a
consequence, some of its global properties too, like for instance the
geodesics equations, which are obtained by the integration of a set of
partial differential equations based upon G.

In \cite{HZ}, pp.149-155, \ a recipe is given for the obtention of G. It
consists on a formula for the norm of tangent vectors (which are systems
derivatives) as a function of the matrices of a state space representation,
needing the solution of a Lyapunov equation; this norm takes
straightforwardly to the metric tensor. Identification algorithms which
require the solution of Lyapunov or Riccati equations are common in the
literature.

Here a formula for the metric tensor is obtained, which gives it directly as
a functional of the system's transfer function. The generality of that
formula enables the use of any representation besides state space (at least
two other ones are known: ARMA and matrix pencil \cite{AP}) and doesn't
require the solution of Lyapunov or Riccati equations. Particularized
formulas for ARMA and state space representations are obtained, using
overlapping parametrizations.

\section{Definitions}

An ARMA(p,q) (Auto-Regressive Moving-Averages) model for a vector time
series $\{y_{t}\}$, y $\in \mathcal{R}^{m}$, t = 0,1,2,3..., is an equation
of the type 
\begin{equation}
A_{0}y_{t}+A_{1}y_{t-1}+...+A_{p}y_{t-p}=B_{0}\varepsilon
_{t}+B_{1}\varepsilon _{t-1}+...+B_{q}\varepsilon _{t-q}
\end{equation}

where the A$_{i}$ and the B$_{i}$ are m$\times $m square matrices and \{$%
\varepsilon _{t}$\} is a gaussian white noise of null mean and covariance R.

In the frequency domain (z-transform), the model becomes 
\begin{equation}
\left( A_{0}+A_{1}z^{-1}+...+A_{p}z^{-p}\right) Y(z)=\left(
B_{0}+B_{1}z^{-1}+...+B_{p}z^{-q}\right) \epsilon (z)
\end{equation}

where 
\begin{equation}
Y(z)=\sum_{i=0}^{\infty }y_{i}z^{-i}
\end{equation}

and 
\begin{equation}
\epsilon (z)=\sum_{i=0}^{\infty }\varepsilon _{i}z^{-i}
\end{equation}

and z is a complex variable (if z is restricted to the unitary circle in the
complex plane, it will have module 1 and, so, will be able to be written as $%
z=e^{\omega \text{i}},$with i=$\sqrt{-1}$ and $\omega \in \mathcal{R}$; in
this case. the z-transform is reduced to the discrete Fourier transform,
which justifies the expression \ ''frequency domain '').

Now, let

\begin{equation}
A(z)=A_{0}+A_{1}z^{-1}+...+A_{p}z^{-p}
\end{equation}

\begin{equation}
B(z)=B_{0}+B_{1}z^{-1}+...+B_{q}z^{-q}
\end{equation}

Then, the frequency domain equation can be written as

\begin{equation}
Y(z)=A(z)^{-1}B(z)\epsilon (z)
\end{equation}

(In the interesting cases, the nondegenerate ones, A(z) is invertible).

The so called transfer function of the model is

\begin{equation}
H(z)=A(z)^{-1}B(z)
\end{equation}

Thus,

\begin{equation}
Y(z)=H(z)\epsilon (z)
\end{equation}

It can be proved that, if \{$\varepsilon _{t}$\}is a gaussian stochastic
process, $\{y_{t}\}$ will also be.

As the definition of H(z) implies that its elements will be fractions whose
numerators and denominators are polynomials in z, that is, rational
functions, the conclusion is that ARMA models represent linear dynamical
systems, that is, systems whose input-output relation is of the type

\begin{equation}
y_{t}=\dsum\limits_{i\text{=}0}^{\infty }H_{i}\varepsilon _{t-i}\bigskip
\end{equation}

where $H(z)=\dsum\limits_{i\text{=}0}^{\infty }H_{i}z^{-i}\ \ \ \ \ $

being the H$_{i}$ p X p square matrices known as the system's Markov
parameters (or weighting sequence). A linear dynamical system is defined by
its sequence $\{H_{i}\}$ of Markov parameters, so that each system s can be
viewed as a point (or vector) in the Hilbert space of these sequences,
defined by:

\textbf{Sum of two systems}: s$_{3}=$s$_{1}+$ s$_{2}$ is the system such
that H$_{i}^{(3)}=$ H$_{i}^{(1)}+$ H$_{i}^{(2)}$ , i = 0,1,2,....

\textbf{Product of a system by a scalar }$\alpha :$ s$_{2}$= $\alpha $s$_{1}$
is the system such that H$_{i}^{(2)}=\alpha $H$_{i}^{(1)}$ , i = 0,1,2,...

\textbf{Internal product between two systems:}

\TEXTsymbol{<}s$_{1},$s$_{2}>$ = tr$\left[ \dsum\limits_{i\text{=}0}^{\infty
}H_{i}^{(1)}(H_{i}^{(2)})^{T}\right] $,

where $T$ represents matrix transposition.

\textbf{Norm of a system:}\TEXTsymbol{\vert}\TEXTsymbol{\vert}s\TEXTsymbol{%
\vert}\TEXTsymbol{\vert} \textbf{= }$\sqrt{<s_{1},s_{2}>}$

\textbf{Distance between two systems:}d(s$_{1}$,s$_{2})$ = \TEXTsymbol{\vert}%
\TEXTsymbol{\vert}s$_{1}-$s$_{2}$\TEXTsymbol{\vert}\TEXTsymbol{\vert}

\smallskip

These definitions in terms of \{H$_{i}\}$ are equivalent to the following
ones in terms of H(z):

\textbf{Sum of two systems}: s$_{3}=$s$_{1}+$ s$_{2}$ is the system such
that H$_{3}(z)=$ H$_{1}(z)+$ H$_{2}(z)$ .

\textbf{Product of a system by a scalar }$\alpha :$ s$_{2}$= $\alpha $s$_{1}$
is the system such that H$_{2}(z)=\alpha $H$_{1}(z)$ , i = 0,1,2,...

\textbf{Internal product between two systems: }

\TEXTsymbol{<}s$_{1},$s$_{2}>$ = $\frac{1}{2\pi \text{i}}\doint_{C}tr\left[
H_{1}(z)H_{2}^{T}(z^{-1})\right] z^{-1}dz$,

where i=$\sqrt{-1}$and C is the unit circle in complex plane.

\textbf{Norm of a system:}\TEXTsymbol{\vert}\TEXTsymbol{\vert}s\TEXTsymbol{%
\vert}\TEXTsymbol{\vert} \textbf{= }$\sqrt{<s_{1},s_{2}>}$

\textbf{Distance between two systems:}d(s$_{1}$,s$_{2})$ = \TEXTsymbol{\vert}%
\TEXTsymbol{\vert}s$_{1}-$s$_{2}$\TEXTsymbol{\vert}\TEXTsymbol{\vert}

The system will be stable if the minimum common multiple of the denominators
of the elements of H(z) (expressed as a rational matrix with all polynomials
in z, not in z$^{-1}$(the B operator of Box \&\ Jenkins)) has all of its
roots inside the unit circle in the complex plane. The system's stability is
a condition that guarantees that to a stationary input process \{$%
\varepsilon _{t}$\} there corresponds a stationary output process $\{y_{t}\}$
, stationarity here understood as time invariance of all the moments of the
stochastic process.

Let, now, the system's Hankel matrix be defined by 
\begin{equation}
\mathcal{H}=\left[ 
\begin{array}{llll}
H_{1} & H_{2} & H_{3} & ... \\ 
H_{2} & H_{3} & H_{4} & ... \\ 
... & ... & ... & ...
\end{array}
\right]
\end{equation}

and search its lines top down, retaining only the ones which are linearly
independent with the preceding ones. Associated to $\mathcal{H}$ there are m
(the number of components of the output) positive integers, the so called
Kronecker indices of the system; the index n$_{i}$ is the number of retained
lines corresponding to the i-th component of the vector time series $%
\{y_{t}\}$ . The so called McMillan degree of the system is the rank of $%
\mathcal{H}$, which will be denoted by n. It is clear, then, that 
\begin{equation}
n=\sum_{i=1}^{m}n_{i}
\end{equation}

\section{Why aren't p and q invariant}

In the univariate case, m=1, A(z) e B(z) are polynomials in z$^{-1}$.
Adopting the convention of using low case letters for scalars, the transfer
function can be written as

\bigskip

\begin{equation}
h(z)=\frac{b(z)}{a(z)}
\end{equation}

\bigskip

If b(z)=p(z)c(z) and a(z)=p(z)d(z), where p(z), c(z) and d(z) are
polynomials in z$^{-1}$, the common factor p(z) can be cancelled, so that
the same system (same transfer function) can be represented by the model

\begin{equation}
c(z)Y(z)=d(z)\epsilon (z)
\end{equation}

If p was the degree of a(z), q of b(z) and r of p(z), the new model will be
an ARMA(p-r,q-r). Thus, cancelling all the common factors of a(z) and b(z),
an ARMA(p*,q*) with p* and q* minimal is obtained. It is obvious that any
increase or reduction in p has, necessarily, to be accompanied by the same
increase or reduction in q, since increases or reductions require adding or
cancelling of common factors in the fraction h(z). The conclusion is that,
in the scalar case, p* and q* are individually invariant, that is, given a
system, it makes sense to refer to \textit{the} order p* of its AR part and 
\textit{the} order q* of its MA part and trying to identify them.

In the vector case, nevertheless, this is not true. In this case, the
transfer function is

\begin{equation}
H(z)=A(z)^{-1}B(z)
\end{equation}

If A(z)=P(z)C(z) and B(z)=P(z)D(z), with P,C e D matrix polynomials, then

\begin{equation}
H(z)=C(z)^{-1}P(z)^{-1}P(z)D(z)=C(z)^{-1}D(z)
\end{equation}

Now, C=P$^{-1}A$ and D=P$^{-1}B.$ If C and D are to be polynomial (so that
CY=D$\epsilon $ be an ARMA model), it is necessary that all elements of A
and B be divisible by detP (since P$^{-1}=\frac{1}{\det P}Cof(P)^{T}).$ Let
d(z) be the maximum common divisor of elements of A e B. Define 
\begin{equation}
P^{-1}=\frac{1}{d(z)}I
\end{equation}

where I is the identity matrix. Then, P=d(z)I is polynomial and the only
polynomial matrices M such that C=M$^{-1}C$ and D=M$^{-1}D$ are polynomial
will be the ones whose determinant is not a polynomial in z$^{-1}$ but a
numerical constant; such polynomial matrices are known as unimodular
matrices. Pairs (C,D) with such a property are said left coprime \cite{KT}.
In the univariate case, after the cancelling of all common polynomial
divisors, there remain only arbitrary numerical constants, which, canceled
between numerator and denominator, don't change their degrees. In the
multivariate case, there remain arbitrary unimodular matrices, which, once
discounted, do alter the AR and MA degrees of the equation. That's why the
following theorem can be stated:

\bigskip

\begin{theorem}
\textbf{\ }In the multivariate case, the minimal p and q are not
individually invariant.
\end{theorem}

\bigskip

\textbf{Proof: }Let s be a system representable by an ARMA(p,q+r) model of
type [A(z),U(z)B(z)] irreducible (that is, with A and UB left coprime), with
U unimodular, detA$_{p}\neq 0$ e detB$_{q}\neq 0$, where 
\begin{equation}
U(z)=U_{0}+U_{1}z^{-1}+...+U_{r}z^{-r}
\end{equation}

\begin{equation}
A(z)=A_{0}+A_{1}z^{-1}+...+A_{p}z^{-p}
\end{equation}

\begin{equation}
B(z)=B_{0}+B_{1}z^{-1}+...+B_{q}z^{-q}
\end{equation}

The restrictions over the determinants imply that the product of A(z) or
B(z) by any polynomial matrix increase their degrees (for instance,
U(z)B(z)=M$_{0}+M_{1}z^{-1}+...M_{q+r}z^{-q+r}$, with M$_{q+r}=U_{r}B_{q}%
\neq O$, since $U_{r}B_{q}=O$ would imply on $%
U_{r}B_{q}B_{q}^{-1}=OB_{q}^{-1}=>U_{r}=O$).

\smallskip The coprimeness of (A,UB) implies on the impossibility of
reducing the degrees of A and UB through the cancelling of a non-unimodular
matrix P such that (A,UB)=(PM,PN). There remain only the unimodular ones.

Now, if U is unimodular, U$^{-1}$ will be polynomial (and, by the way,
unimodular too). So, U$^{-1}A$ will also be polynomial.

The conclusion is that (U$^{-1}A,U^{-1}UB)=(U^{-1}A,B)$ will be an
ARMA(p+k,q) model, where k=degree(U$^{-1})$, irreducible, which will
represent the same system. So, there are systems which have two irreducible
ARMA representations with different degrees.

\textbf{Q.E.D.}

Consider, for instance, a system whose transfer function H(z) is itself a
unimodular matrix. Then, this system has an ARMA(0,q) representation of the
type 
\begin{equation}
Y(z)=H(z)\epsilon (z)
\end{equation}

with 
\begin{equation}
H(z)=B_{0}+B_{1}z^{-1}+...+B_{r}z^{-r}
\end{equation}

But, in this case, H(z)$^{-1}$ is also polynomial, since 
\begin{equation}
H(z)^{-1}=\frac{1}{\det H(z)}Cof[H(z)]^{T}
\end{equation}

since det(Hz) is a number and its cofactor matrix is polynomial.

Then, multiplying the ARMA(0,q) equation above by H(z)$^{-1}$, an ARMA(q,0)
equation is obtained:

\begin{equation}
H(z)^{-1}Y(z)=H(z)^{-1}H(z)\epsilon (z)=\epsilon (z)
\end{equation}

Thus, there are two pairs, (0,q) and (q,0), of minimal orders (because nor p
nor q can be smaller then zero) corresponding to the same system. Hannan and
Deistler (\cite{HD}, pg. 77), although don't call attention to this
phenomenon, exhibit,\textit{\ en passant}, an example with q=1.

\begin{remark}
\textbf{\ } Although the systems built in the proof have a minimum p* and a
minimum q*, this is \ not useful, since they don't have an ARMA(p*,q*)
representation. In the case (0,q)-(q,0) exhibited, for instance, the system
doesn't have an ARMA(0,0) representation (save for very special cases),
which would represent a white noise. What would be useful for model building
is the joint minimality of p and q, which would imply on the invariance of
the structure of the minimum \ ARMA model.
\end{remark}

\section{The ARMA(p,p)-irreducible class}

Given the impossibility of representing any system by an ARMA(p,q) model
with p and q jointly minimal, something that can be done is to represent it
by an ARMA(p,p), that is, with p=q, such that p be the least possible
integer. It can be proved (see \cite{HZ}, pg. 38) that this minimum value of
p is equal to the largest of the Kronecker indices of the system.

The set of all systems representable by an ARMA(p,p) model with $%
A(z)=A_{0}+A_{1}z^{-1}+...+A_{p}z^{-p}$, with p equal to their largest
Kronecker \ index- that is, with p=max\{n$_{i}\}$ - will be called \emph{\
ARMA(p,p)-irreducible class} (the word irreducible here relates to the
impossibility of reducing the value of p). It is a subset of the Hilbert
space of linear systems. Some subsets of it are made up of systems with more
than one ARMA(p,p)-irreducible model; so, the ARMA(p,p)-irreducible
parametrization is not identifiable in the strict sense (although it is in
Vajda's sense). (There are here, as always, two sets: the set of ARMA(p,p)
models and the set of systems with max\{n$_{i}\}=p$, which are the image of
those models in Hilbert's space. The parametrization which maps a set on the
other one is said to be identifiable if it is biunivocal).

The problem posed here is if, analogously to the set of systems with a
common $\ \sum_{i=1}^{m}n_{i}$ , the set of systems with a common max\{n$%
_{i}\}$ is a differentiable manifold; for, if that was the case, it would be
possible to cover it with a set of charts, thus obtaining an overlapping
parametrization which would be the most natural for the ARMA representation
(in the state space representation, the natural thing is to treat with the
set of systems with a common $\sum_{i=1}^{m}n_{i}$ , which is the minimal
dimension of the state space). To start the analysis, consider, firstly, an
ARMA canonical form for systems with Kronecker indices n$_{i},$ i=1...m. It
can be obtained by the following procedure:

Aline the components of the predictor y(t/t-1) (y(t) conditioned to the time
series until time t-1), $y_{1},y_{2},...,y_{m}$, from t-p to p. Next, adopt
the procedure indicated for the example below, in which m=3, $%
n_{1}=2,n_{2}=2,n_{3}=2:$

\smallskip 
\begin{equation}
\left[ 
\begin{array}{lll}
y_{1} & y_{2} & y_{3} \\ 
\ast & * & * \\ 
&  & 
\end{array}
\right] 
\begin{array}{lll}
y_{1} & y_{2} & y_{3} \\ 
\ast & o & * \\ 
\ast & * & *
\end{array}
\left[ 
\begin{array}{lll}
y_{1} & y_{2} & y_{3} \\ 
o &  & o \\ 
\ast &  & 
\end{array}
\right]
\end{equation}

\smallskip

where the leftmost matrix represents time t-2, the central one, t-1, and the
right one, t. In the second line,the \ ''x '' indicate the components which
enter in the selection and the ''o '' indicate the first time, in the search
from left to right, that a component revealed to be linearly dependent with
the precedent ones. The third line is the second one dislocated to the left,
so as to have all L.D. components in the same column (the third one).

Now consider the third column. The first component (indicated in the second
line) is L.D. with all components two time units behind and with the first \
and the third ones one time unit behind; so, in general, there will be
coefficients corresponding to these components when writing down the model:

The second component depends of the first one in the same time and of all
one time unit behind, that is,

$y_{2}(t/t-1)=a_{2}y_{1}(t/t-1)+b_{2}y_{1}(t-1/t-2)$

$+c_{2}y_{2}(t-1/t-2)+d_{2}y_{3}(t-1/t-1)$

The third one depends on the first and the third ones one time unit behind
and of all two time units behind, that is,

$y_{3}(t/t-1)=a_{3}y_{1}(t-1/t-2)+b_{3}y_{3}(t-1/t-2)$

$+c_{3}y_{1}(t-2/t-3)++d_{3}y_{2}(t-2/t-3)+e_{3}y_{3}(t-2/t-3)$

\smallskip As to the MA part, all of its matrices are full, with the only
restriction that, in each line of the equation, the degree of the MA part
cannot be greater than the degree of the AR part.

Thus, the canonical structure becomes:

$\left[ 
\begin{array}{lll}
1 & 0 & 0 \\ 
x & 1 & 0 \\ 
0 & 0 & 1
\end{array}
\right] y_{t}+\left[ 
\begin{array}{lll}
x & 0 & x \\ 
x & x & x \\ 
x & 0 & x
\end{array}
\right] y_{t-1}+\left[ 
\begin{array}{lll}
x & x & x \\ 
0 & 0 & 0 \\ 
x & x & x
\end{array}
\right] y_{t-2}=$

$=\varepsilon _{t}+\left[ 
\begin{array}{lll}
x & x & x \\ 
x & x & x \\ 
x & x & x
\end{array}
\right] \varepsilon _{t-1}+\left[ 
\begin{array}{lll}
x & x & x \\ 
0 & 0 & 0 \\ 
x & x & x
\end{array}
\right] \varepsilon _{t-2}$

Formally, what one has in the general case is (see \cite{GD}, here modified
to the stochastic case):

\begin{equation}
y_{i}(t+n_{i})=\sum_{j=1}^{m}\sum_{k=1}^{n_{ij}}a_{ijk}y_{j}(t+k-1)+%
\sum_{j=1}^{m}\sum_{k=1}^{n_{i}}b_{jik}\varepsilon _{j}(t+k-1)+\varepsilon
_{i}(t+n_{i})
\end{equation}

where n$_{ij}=\left( 
\begin{array}{l}
n_{i}\text{, for i=j} \\ 
\min \{n_{i}+1,n_{j}\}\text{, for i\TEXTsymbol{>}j } \\ 
\min \{n_{i},n_{j}\}\text{, for i\TEXTsymbol{<}j }
\end{array}
\right) $

shifting, next, each equation in time, so that in the left member always
appear $y_{i}(t).$

The important to consider here are the following three facts:

1) detA$_{0}=1$ always, so that one can multiply all the vector equation by A%
$_{0}^{-1}$, to obtain a monic model, that is, with A$_{0}=$ I (mXm
identity), without loss of generality.

2) Only when all Kronecker indices are equal is that all matrices - except
for A$_{o}$, which will be the identity - are full. This is, thus, the case
with the greatest number of free parameters: 2m$^{2}$p; the dimension of its
image M in the Hilbert space of systems will be, thus, also 2m$^{2}$p. It is
known that M is a differentiable manyfold, since each set of systems with
the same m (number of components of y) and same set of Kronecker indices is
the image of one of the maps of the chart which defines a differentiable
manyfold of dimension 2mn \cite{CK}, where $n=\sum_{i=1}^{m}n_{i}$ (when all
indices are equal, one has $p=n_{1}=...=n_{m}$; thus n=mp and so 2m$%
^{2}p=2mn $). This is, consequently, the generic case of this
parametrization.

3) In all the other cases, the augmented matrix $[A_{p},B_{p}]$ will not be
full rank (its product by an invertible matrix - A$_{0}^{-1}$, for instance
- thus, also not).

\begin{definition}
The ARMA(p,p)-irreducible parametrization will be defined by
\end{definition}

\begin{equation}
y_{t}+A_{1}y_{t-1}+...+A_{p}y_{t-p}=\varepsilon _{t}+B_{1}\varepsilon
_{t-1}+...+B_{p}\varepsilon _{t-p}
\end{equation}

with all matrices, in principle, full. The word irreducible denotes the
impossibility of reducing p. Thus, such parametrization includes all systems
with max\{n$_{i})=p$ and only them.

From the properties of the canonical forms here exhibited, follows the
conclusion that the image of the set of ARMA(p,p)-irreducible models in the
space of systems is the union of a differentiable manyfold of dimension 2pm$%
^{2}$ (corresponding to all systems with $n_{1}=...=n_{m}=p)$ with sets of
lower dimensions (corresponding to systems with some n$_{i}\neq p)$. For
easier references, it is convenient to state the following

\begin{theorem}
The generic sub-class of the ARMA(p,p)-irreducible

parametrization is a differentiable manyfold of dimension 2pm$^{2}$, where m
is the number of components of $\left\{ y_{t}\right\} .$
\end{theorem}

\textbf{Proof: }The generic sub-class of this parametrization is the one in
which all matrices are full and rank$\left[ A_{p},B_{q}\right] =m.\,$But in
this case the ARMA model is under the canonical form of systems with $%
n_{1}=...=n_{m}=p$, which, as is known \cite{CK}, is a differentiable
manifold of dimension 2pm$^{2}$.

\textbf{Q.E.D.}

The problem here considered is of knowing if that union of systems sets,
that is, the complete image, constitutes a differentiable manifold (in this
case, its dimension would be 2pm$^{2}$)$.$ This will occur if the points
(systems) of the non-generic sets were regular under any coordinates system,
which would require that the set of tangent vectors at each point spanned a
space of dimension exactly equal to 2pm$^{2}.$ Unfortunately, the following
theorem shows that this is not the case:

\begin{theorem}
\textbf{\ }The set of systems representable by the

ARMA(p,p)-irreducible parametrization does not constitute a differentiable
manyfold.
\end{theorem}

\smallskip

\textbf{Proof:} Let s be a system representable by this parametrization with
some n$_{i}\neq p$. Let M(z$)=I+Mz^{-1}$ with $M[A_{p},B_{p}]=O$ (null
matrix). For instance, in the case above illustrated$%
(m=3,n_{1}=2,n_{2}=2,n_{3}=2),$

\begin{equation}
M[A_{p},B_{p}]=\left[ 
\begin{array}{lll}
0 & a & 0 \\ 
0 & b & 0 \\ 
0 & c & 0
\end{array}
\right] \left[ 
\begin{array}{llllll}
x & x & x & x & x & x \\ 
0 & 0 & 0 & 0 & 0 & 0 \\ 
x & x & x & x & x & x
\end{array}
\right] =\left[ 
\begin{array}{lll}
0 & 0 & 0 \\ 
0 & 0 & 0 \\ 
0 & 0 & 0
\end{array}
\right]
\end{equation}

and

\begin{equation}
M(z)=\left[ 
\begin{array}{lll}
1 & 0 & 0 \\ 
0 & 1 & 0 \\ 
0 & 0 & 1
\end{array}
\right] +\left[ 
\begin{array}{lll}
0 & a & 0 \\ 
0 & b & 0 \\ 
0 & c & 0
\end{array}
\right] z^{-1}=\left[ 
\begin{array}{lll}
1 & az^{-1} & 0 \\ 
0 & 1+bz^{-1} & 0 \\ 
0 & cz^{-1} & 1
\end{array}
\right]
\end{equation}

As in the systems with some n$_{i}\neq p$ the extended matrix $[A_{p},B_{p}]$
is not full rank, there always is M with the properties above. In this case,
the pair $\left[ M(z)A(z),M(z)B(z)\right] $ will also represent an
ARMA(p,p)-irreducible model with the same p, for the same system s, whatever
be the values of a,b and c.

Let, now, H$_{\theta _{i}}(z)$, i=1...2pm$^{2}$, be the tangent vectors at a
system representable by an ARMA(p,p)-irreducible parametrization, where the
free parameters were denoted by $\{\theta _{i},i=1,2,...\}$. Since H(z) can
be factored in

$H(z)=A(z)^{-1}B(z)$

being A(z) and B(z) the polynomial matrices

\begin{equation}
A(z)=I+A_{1}z^{-1}+...+A_{p}z^{-p}
\end{equation}

\begin{equation}
B(z)=I+B_{1}z^{-1}+...+B_{p}z^{-p}
\end{equation}

it follows that

\begin{equation}
H_{\theta _{i}}(z)=A^{-1}(z)B_{\theta _{i}}(z)-A^{-1}(z)A_{\theta
_{i}}(z)A^{-1}(z)B(z)
\end{equation}

Or,

\begin{equation}
H_{\theta _{i}}(z)=A^{-1}(z)\left( B_{\theta _{i}}(z)-A_{\theta
_{i}}(z)H(z)\right)  \label{deriv}
\end{equation}

Now, $B_{\theta _{i}}(z)$ and $A_{\theta _{i}}(z)$ are constant matrices
(consider, for instance, the case m=2, p=2):

$y_{t}+\left[ 
\begin{array}{ll}
\theta _{1} & \theta _{2} \\ 
\theta _{3} & \theta _{4}
\end{array}
\right] y_{t-1}+\left[ 
\begin{array}{ll}
\theta _{5} & \theta _{6} \\ 
\theta & \theta _{8}
\end{array}
\right] y_{t-2}=$

$=\varepsilon _{t}+\left[ 
\begin{array}{ll}
\theta _{9} & \theta _{10} \\ 
\theta _{11} & \theta _{12}
\end{array}
\right] \varepsilon _{t-1}+\left[ 
\begin{array}{ll}
\theta _{13} & \theta _{14} \\ 
\theta _{15} & \theta _{16}
\end{array}
\right] \varepsilon _{t-2}$

Then, one has, for instance,

$A_{\theta _{3}}(z)$ =$\left[ 
\begin{array}{ll}
0 & 0 \\ 
1 & 0
\end{array}
\right] z^{-1}$

$B_{\theta _{13}}(z)=\left[ 
\begin{array}{ll}
1 & 0 \\ 
0 & 0
\end{array}
\right] z^{-2}$ )

The matrix H(z) is unique and can be written in terms only of the free
parameters of the canonical form of s.

However, if s is of the kind here considered (that is, such that some n$%
_{i}\neq p$), it will have an infinite number of representations $\left[ 
\mathcal{A}\text{,}\mathcal{B}\right] $, all of them being
ARMA(p,p)-irreducible with the same p, as already claimed, of the type

$\mathcal{A}(z)=M(z)A(Z)$

$\mathcal{B}(z)=M(z)B(z)$

where A(z) and B(z) are the canonical matrices.

This doesn't change H(z), since

$\mathcal{H}(z)=\mathcal{A}(z)^{-1}\mathcal{B}(z)=[M(z)A(z)]^{-1}[M(z)B(z)]$

$=A(z)^{-1}M^{-1}(z)M(z)B(z)=A(z)^{-1}B(z)=H(z)$

But it does change A(z) in (\ref{deriv}) and, so, H$_{\theta _{i}}$,
i=1...2pm$^{2}.$ In particular, it is always possible to define M(z)
non-unimodular, maintaining the properties above (the case exhibited is an
example: it suffices that b$\neq 0$), such that M$^{-1}(z)$ , when
right-multiplied by A$^{-1}(z)$ in (\ref{deriv}) introduce a free parameter
(in this instance, b) in the denominator of H$_{\theta _{i}}.$ Then, any
change in the value of this parameter will result in a tangent vector (a
system) which is L.I. with the remaining ones. Indeed, calling $\gamma $ =($%
\gamma _{1},\gamma _{2},...,\gamma _{r})$ the vector of free parameters of
M(z) which appear in its determinant, one concludes that the space spanned
by the tangent vectors at s is of infinite dimension (thus, different from
2pm$^{2})$, that is,

dim\{$H_{\theta _{i}}(z)/\gamma \in R^{r}\}=\infty $

The conclusion is that points s of the kind here considered (that is, those
ones not belonging to the generic sub-class) cannot be made regular under
any coordinates system. Thus, the union of all classes of systems
represented by the ARMA(p,p)-irreducible parametrization is not a
differentiable manifold.

\textbf{Q.E.D.}

\begin{remark}
:The systems belonging to the non-generic case are in a position
qualitatively similar to the vertex of a two-sided cone: it is not a regular
point of the surface, but tangent vectors are perfectly definable on it, the
anomaly being that they span of the R$^{3}$ instead of just a plane
\end{remark}

\section{Riemannian Metric Tensor}

Consider a parametrization (D,P) which defines a class of dynamical systems,
so that the elements h$_{ij}(z)$ of the transfer function matrix H(z) which
defines a system be expressed as a function of $\theta =(\theta _{1},\theta
_{2},...,\theta _{d})$, that is, of the local coordinates of one of the maps
of the differentiable manifold defined by (D,P), whose dimension will be
denoted by ''d ''.

\begin{theorem}
The element (i,j) of the Riemannian metric tensor G of a dynamical linear
system is given, in terms of its transfer function matrix H(z) and an
adequate parametrization that expresses it as a function of a finite vector $%
\theta \in R^{d}$, by:
\end{theorem}

\medskip

\[
g_{ij}=\frac{1}{2\pi \text{i}}\doint_{C}tr\left[ \frac{\partial H(z)}{%
\partial \theta _{i}}\frac{\partial H^{T}(z^{-1})}{\partial \theta _{j}}%
\right] z^{-1}dz 
\]

\medskip

\medskip

where i=$\sqrt{-1}$ , C is the unit circle centered on the origin of the
complex plane and tr stands for trace.

\medskip

\medskip

\textbf{Proof:}

\medskip

Consider a system $s$ belonging to the manifold defined by a

parametrization (D,P). Specifically, suppose that, in one of it's maps, $s$
be represented by the vector of local coordinates $\theta =($ $\theta _{1},$ 
$\theta _{2},...,$ $\theta _{d})$. Let $\theta \left( v_{1}\right) $ and $%
\theta \left( v_{2}\right) $ be two curves of M passing by $s.$ Let also $%
s_{1}^{\prime }$ and $s_{2}^{\prime }$ be the two derivatives of $s$ along
the two curves, respectively.

The Riemannian metric tensor G at $s$ is, as already defined, the matrix
such that

$\theta _{1}^{\prime \text{ T}}G$ $\theta _{2}^{\prime }=$ g$(s_{1}^{\prime
},s_{2}^{\prime })=<s_{1}^{\prime },s_{2}^{\prime }>=$ tr$\left[
\dsum\limits_{i=0}^{\infty }H_{i}^{1}(H_{i}^{2})^{T}\right] $

where H$_{i}^{1}$ and H$_{i}^{2}$ , $i=0,1,2,3,...,$ are the Markov
parameters of the two tangent systems.

By Percival's formula (see \cite{OP}, for the scalar case; the
generalization used here is easily obtained), one has

\medskip

\medskip

tr$\left[ \dsum\limits_{i=0}^{\infty }H_{i}^{1}(H_{i}^{2})^{T}\right] =$ $%
\frac{1}{2\pi \text{i}}\doint_{C}tr\left[ H_{1}(z)H_{2}^{T}(z^{-1})\right]
z^{-1}dz$

\medskip

\medskip

\medskip

where H$_{1}(z)$ and H$_{2}\left( z\right) $ are the transfer functions of
the two tangent systems. Now, these transfer functions, because they refer
to tangent systems along those curves (and recalling that stability of the
systems implies convergence of the infinite series involved, thus justifying
the commutation between summations and derivatives), are given by

\medskip

$H_{k}\left( z\right) =\dsum\limits_{i=0}^{\infty }\frac{dH_{i}}{dv_{i}}%
z^{-i}=\frac{d}{dv_{i}}\dsum\limits_{i=0}^{\infty }H_{i}z^{-i}=\frac{dH(z)}{%
dv_{k}}=\dsum\limits_{i\text{=}1}^{d}\frac{\partial H(z)}{\partial \theta
_{i}}\frac{d\theta _{i}}{dv_{k}}=$

$=\dsum\limits_{i\text{=}1}^{d}\frac{\partial H(z)}{\partial \theta _{i}}%
\theta _{i}^{\prime (k)}$ , k=1,2

\medskip

\medskip

where H(z) is the transfer function of s and $\frac{d\theta _{i}}{dv_{k}}$
was abbreviated to $\theta _{i}^{\prime (k)}.$

Substituting in the last integral, there results

$\medskip $

$\theta _{1}^{\prime \text{ T}}G$ $\theta _{2}^{\prime }=$ $\frac{1}{2\pi 
\text{i}}\doint_{C}tr\left[ \left( \dsum\limits_{i\text{=}1}^{d}\frac{%
\partial H(z)}{\partial \theta _{i}}\theta _{i}^{\prime (1)}\right) \left(
\dsum\limits_{i\text{=}1}^{d}\frac{\partial H(z^{-1})}{\partial \theta _{i}}%
\theta _{i}^{\prime (2)}\right) ^{T}\right] z^{-1}dz$

\medskip

where $\;\theta _{i}^{\prime (1)}$ $=$ $\left[ \theta _{1}^{\prime \text{ }}%
\right] _{i}$ \ \ , i=1,2,...,d

and $\;\;\;\theta _{i}^{\prime (2)}=\left[ \theta _{2}^{\prime \text{ }}%
\right] _{i}$ \ \ , i=1,2,...,d

\medskip which specifies the components of vectors $\theta _{1}^{\prime 
\text{ }}$ and $\theta _{2}^{\prime \text{ }}.$

\medskip

The element g$_{ij}$ of tensor G is obtained by

g$_{ij}=e_{i}^{T}Ge_{j}$

where $e_{i}=\left[ 0,...,0,1,0,...,0\right] $ with the ''1 '' in the i-th
position.

Substitute, in the last integral, $\theta _{1}^{\prime \text{ }}$ by $e_{i}$
and $\theta _{2}^{\prime \text{ }}$ by $e_{j}.$

The summations become:

\medskip

$\dsum\limits_{n=1}^{d}\frac{\partial H(z)}{\partial \theta _{n}}\theta
_{n}^{\prime (1)}=\frac{\partial H(z)}{\partial \theta i}$ \ \ \ \ \ \ \ and
\ \ \ \ \ \ $\dsum\limits_{n=1}^{d}\frac{\partial H(z^{-1})}{\partial \theta
_{n}}\theta _{n}^{\prime (2)}=\frac{\partial H(z^{-1})}{\partial \theta j}$

\medskip

and there results

\medskip

$\fbox{g$_{ij}=\frac{1}{2\pi \text{i}}\doint_{C}tr\left[ \frac{\partial H(z)%
}{\partial \theta _{i}}\frac{\partial H^{T}(z^{-1})}{\partial \theta _{j}}%
\right] z^{-1}dz$}$

\medskip \textbf{Q.E.D.}

\begin{example}
For the scalar ARMA(1,1) parametrization,
\end{example}

$y_{t}+ay_{t-1}=u_{t}+bu_{t-1}$, so that $\theta =(\theta _{1},\theta
_{2})=(a,b),$ with $a,b\in (-1,1)$

and $a\neq b.$ The transfer function in this case is scalar, given by

\medskip

$H(z)=\frac{z+b}{z+a}=\frac{z+\theta _{1}}{z+\theta _{2}}$

\medskip

so that

\medskip

g$_{ij}=\frac{1}{2\pi \text{i}}\doint_{C}tr\left[ \frac{\partial H(z)}{%
\partial \theta _{i}}\frac{\partial H^{T}(z^{-1})}{\partial \theta _{j}}%
\right] z^{-1}dz=$

$=$.$\frac{1}{2\pi \text{i}}\doint_{C}\left( \frac{\partial }{\partial
\theta _{i}}\frac{z+\theta _{1}}{z+\theta _{2}}\right) \left( \frac{\partial 
}{\partial \theta _{j}}\frac{z^{-1}+\theta _{1}}{z^{-1}+\theta _{2}}\right)
z^{-1}dz$

$i,j=1,2$

\medskip Thus,

\medskip

g$_{11}=\frac{1}{2\pi \text{i}}\doint_{C}\frac{z+b}{(z+a)^{2}}\frac{z^{-1}+b%
}{(z^{-1}+a)^{2}}z^{-1}dz$ $=\frac{4ab-(b^{2}+1)(a^{2}+1)}{(a^{2}-1)^{3}}$

\medskip

g$_{12}=$g$_{21}=\frac{1}{2\pi \text{i}}\doint_{C}\frac{z+b}{(z+a)^{2}}\frac{%
1}{z^{-1}+a}z^{-1}dz=\frac{ab-1}{(1-a^{2})^{2}}$

\medskip

g$_{22}=\frac{1}{2\pi \text{i}}\doint_{C}\frac{1}{z+a}\frac{1}{z^{-1}+a}%
z^{-1}dz=\frac{1}{1-a^{2}}$

\medskip

As a result, the Riemannian metric tensor in this case is

\medskip

\medskip

G = $\left[ 
\begin{array}{ll}
\frac{4ab-(b^{2}+1)(a^{2}+1)}{(a^{2}-1)^{3}} & \frac{ab-1}{(1-a^{2})^{2}} \\ 
\frac{ab-1}{(1-a^{2})^{2}} & \frac{1}{1-a^{2}}
\end{array}
\right] $

\medskip

\medskip which agrees with the corresponding state-space case given in [1,
pg. 222].

\bigskip

\subsection{Overlapping Parametrizations}

\medskip

Consider the case in which H$_{0}=$ I (identity matrix) and y$_{t}$ and u$%
_{t}\in R^{m}$, that is, r=m (number of inputs equal to number of outputs).
There is no loss of generality, since all required is static redefinition of
the inputs (or outputs) and the introduction of artificial dummy inputs (or
outputs).

Consider, now, the differentiable manifold S$_{n}$ of all systems of this
kind with McMillan degree n fixed.\ For m$\neq 1,$ it is not possible to
cover it with a unique map, so that the global chart is made up of a set of
maps, each one of them characterized by a set of m natural numbers n$_{i},$
i=1,2,...,m, with $\sum\limits_{i=1}^{m}n_{i}=n.$ Call M$_{n\text{; }%
n_{1},n_{2},...,n_{m}}$ the corresponding map. Then, the coordinates of a
system s described by this map are defined by the following procedure \cite
{GW}:

Given the system's Hankel matrix

\bigskip

\bigskip

$\mathcal{H}$ = $\left[ 
\begin{array}{llll}
H_{1} & H_{2} & H_{3} & ... \\ 
H_{2} & H_{3} & H_{4} & ... \\ 
... & ... & ... & ...
\end{array}
\right] $

\medskip

\bigskip

let H$^{i}$ be its i-th block of m lines (for instance, H$%
^{2}=[H_{2}H_{3}H_{4}...]$) and let $h_{1i},h_{2i},...,h_{mi}$ be the lines
(of infinite size) of H$^{i}$. The fact that a system can be represented by
the map M$_{n\text{; }n_{1},n_{2},...,n_{m}}$ means that the lines

$h_{11},...,h_{1n_{1}};h_{21,...,}h_{2n_{2}};h_{m1},...h_{mn_{m}}$
constitute a base for the space of lines of $\mathcal{H}$. Thus, the lines $%
h_{1(n_{1}+1)},...,h_{m(n_{m}+1)}$ can be written as linear combinations of
them:

\begin{equation}
h_{i(n_{i}+1)}=\sum\limits_{j=1}^{m}\sum\limits_{k=1}^{n_{j}}\alpha
_{ijk}h_{jk}
\end{equation}
\ \ \ \ \ \ \ \ \ i=1,...,m.

Now, call $h_{ij}(k)$ the k-th element of line h$_{ij}.$ The 2mn numbers \{$%
\alpha _{ijk},k=1,...,n_{j};i,j=1,...,m,$ and $%
h_{ij}(k),i=1,...,m;j=1,...n_{i};k=1,...,p$\} are the system's coordinates
according to map M$_{n\text{; }n_{1},n_{2},...,n_{m}},$ that is, the
components $\theta _{i}$ $\left( i=1,2,...,2nm\right) $ of the
coordinate-vector $\theta $ which appear in the tensor's formula given by
theorem 1. Let the first nm components of $\theta $ be the $\{\alpha
_{ijk}\} $ and the last nm ones be the $\{h_{ij}(k)\}.$

\bigskip

Now, define the matrix K with n rows and m columns as in \cite{GW}.

Each element of this matrix is one of the coordinates $h_{ij}(k).$

\bigskip

\bigskip

\subsection{Tensor for the ARMA representation}

\bigskip

The autoregressive moving averages representation (ARMA)

\bigskip $%
A_{0}y_{t}+A_{1}y_{t-1}+...+A_{p}y_{t-p}=B_{0}u_{t}+B_{1}u_{t-1}+...+B_{p}u_{t-p} 
$

where y$_{t}$,u$_{t}$ $\in R^{m}$ \ , \ A$_{i}$ and B$_{i}$ are $m\times m$
matrices and whose

transfer function is H(z)=A$^{-1}\left( z\right) B(z)$

where \ \ \ \ \ \ \ \ \ A(z)= A$_{0}z^{p}+A_{1}z^{p-1}+...+A_{p}\mathstrut
\strut \bigskip $

and \ \ \ \ \ \ \ \ \ \ \ \ \ \ \ B(z)= B$_{0}z^{p}+B_{1}z^{p-1}+...+B_{p}\ $

has an overlapping parametrization as defined in \cite{GW} in which

p=max\{n$_{i}\}$

a$_{ii}(z)=z^{n_{i}}-\alpha _{iin_{i}}z^{n_{i}-1}-...-\alpha _{ii1}$

a$_{ij}(z)=-\alpha _{ijn_{j}}z^{n_{j}-1}-\alpha
_{ijn_{j}-1}z^{n_{j}-2}-\alpha _{ij1}$, \ \ \ \ i$\neq j$

\bigskip

with the $\alpha _{ijk}$ already defined, and

\bigskip

B(z) = A(z) + M(z)K

\bigskip

where M(z) is a polynomial matrix whose entries are

pseudo-derivatives (in relation to z) of the entries of A(z) (if f(z)=z$%
^{m}+a_{1}z^{m-1}+...+a_{m}$ , then f$^{(k)}\left( z\right)
=z^{m-k}+a_{1}z^{m-k-1}+...+a_{m-k}$ is its pseudo-derivative of order k).

With the definition of vector $\theta $ given above, notice that matrices
A(z) and M(z) are functions only of its first nm components, while matrix K
is function only of its nm last ones. Thus,the following developments can be
made:

\bigskip

H(z) = A$^{-1}$(z)B(z)

\bigskip

$\frac{\partial }{\partial \theta }H(z)=A^{-1}(z)\left[ \frac{\partial }{%
\partial \theta }B(z)-\frac{\partial A(z)}{\partial \theta }H(z)\right] $

\bigskip

$\frac{\partial }{\partial \theta }B(z)=\frac{\partial }{\partial \theta }%
A(z)+\frac{\partial }{\partial \theta }M(z)K+M(z)\frac{\partial }{\partial
\theta }K$

\bigskip \bigskip with

$\frac{\partial }{\partial \theta _{i}}A(z)=\frac{\partial }{\partial \theta
_{i}}M(z)=0$ \ \ \ for \ \.{i} \TEXTsymbol{>} nm \ \ \ \ 

\medskip

$\frac{\partial }{\partial \theta _{i}}K=0$ \ \ \ \ \ \ \ \ \ \ \ \ \ \ \ \
for \ i \TEXTsymbol{<} nm+1

\medskip

\bigskip

So, \ for i \TEXTsymbol{>} nm, there results

\bigskip

\bigskip

$\frac{\partial }{\partial \theta _{i}}H(z)=A^{-1}(z)\left[ \frac{\partial }{%
\partial \theta _{i}}B(z)-\frac{\partial A(z)}{\partial \theta _{i}}H(z)%
\right] =A^{-1}(z)\frac{\partial }{\partial \theta _{i}}B(z)=$

\medskip

$=A^{-1}(z)M(z)\frac{\partial }{\partial \theta _{i}}K$

\bigskip

\bigskip

And, for i \TEXTsymbol{<} nm+1,

\bigskip

$\frac{\partial }{\partial \theta _{i}}H(z)=A^{-1}(z)\left[ \frac{\partial }{%
\partial \theta _{i}}A(z)+\frac{\partial M(z)}{\partial \theta _{i}}K-\frac{%
\partial A(z)}{\partial \theta _{i}}H(z)\right] =$

\bigskip

$=A^{-1}(z)\left( \frac{\partial M(z)}{\partial \theta _{i}}K+\frac{\partial
A(z)}{\partial \theta _{i}}\left[ I-H(z)\right] \right) $

\bigskip

\medskip The results above demand a separation in three cases, each one
giving rise to a corresponding formula for the metric tensor.

Let I=\{1,2,...,nm\} and J=\{nm+1,nm+2,...,2nm\} be index sets.

\bigskip

\bigskip Case 1: \ \ \ \ \ i,j$\in $J

\bigskip g$_{ij}=\frac{1}{2\pi \text{i}}\doint_{C}tr\left[ A^{-1}(z)M(z)%
\frac{\partial K}{\partial \theta _{i}}\frac{\partial K^{T}}{\partial \theta
_{j}}M^{T}(z^{-1})A^{-T}(z^{-1})\right] z^{-1}dz$

\bigskip Case 2: \ \ \ \ \ i,j$\in $I

\bigskip g$_{ij}=\frac{1}{2\pi \text{i}}\doint_{C}tr[A^{-1}(z)\left( \frac{%
\partial M(z)}{\partial \theta _{i}}K+\frac{\partial A(z)}{\partial \theta
_{i}}\left[ I-H(z)\right] \right) .$

.$\left( K^{T}\frac{\partial M^{T}(z^{-1})}{\partial \theta _{j}}+\left[
I-H^{T}(z^{-1})\right] \frac{\partial A^{T}(z^{-1})}{\partial \theta _{j}}%
\right) A^{-T}(z^{-1})]z^{-1}dz$

\bigskip

Case 3: \ \ \ \ \ i$\in $I \ \ and \ \ j$\in $J \ \ \ 

\bigskip

\bigskip g$_{ij}=\frac{1}{2\pi \text{i}}\doint_{C}tr$

$\left[ A^{-1}(z)\left( \frac{\partial M(z)}{\partial \theta _{i}}K+\frac{%
\partial A(z)}{\partial \theta _{i}}\left[ I-H(z)\right] \right) \frac{%
\partial K^{T}}{\partial \theta _{j}}M^{T}(z^{-1})A^{-T}(z^{-1})\right]
z^{-1}dz$

\bigskip

\subsection{Tensor for the state-space representation}

\bigskip

The state space representation

\bigskip

$x_{t+1}=Ax_{t}+Bu_{t}$

$y_{t}=Cx_{t}+u_{t}$

\bigskip

where $x_{t}\in R^{n}$ and $u_{t},y_{t}\in R^{m}$

has the overlapping parametrization defined in \cite{GW}, where

C is a matrix of zeros and ones, A is a sparse matrix in which the only
non-nul variable entries are the $\alpha _{ijk}$ already defined and $B=K$
(as defined above).

Thus,the following calculations will provide explicit formulas for the
tensor:

$\bigskip $

$H(z)=C(zI-A)^{-1}B+I$

\bigskip

and $\frac{\partial }{\partial \theta }$ \bigskip H(z)=$\frac{\partial }{%
\partial \theta }$C(zI-A)$^{-1}B$

and, for i \TEXTsymbol{>} nm,

\bigskip $\frac{\partial }{\partial \theta _{i}}$ \bigskip H(z)=$\frac{%
\partial }{\partial \theta _{i}}$C(zI-A)$^{-1}B=$C(zI-A)$^{-1}\frac{\partial 
}{\partial \theta _{i}}B$

\bigskip whilst, for i \TEXTsymbol{<} nm+1,

$\frac{\partial }{\partial \theta _{i}}$ \bigskip H(z)=$\frac{\partial }{%
\partial \theta _{i}}$C(zI-A)$^{-1}B=$C$\left[ \frac{\partial }{\partial
\theta _{i}}(zI-A)^{-1}\right] B=$

$=C(zI-A)^{-1}\frac{\partial A}{\partial \theta _{i}}(zI-A)^{-1}$

\bigskip Again, this generates three formulas for the tensor:

\bigskip

Case 1: \ \ \ \ i,j$\in $J \ \ \ \ \ (see definitions of I and J in the
preceding section)

\bigskip

$g_{ij}=\frac{1}{2\pi \text{i}}\doint_{C}tr\left[ C(zI-A)^{-1}\frac{\partial
B}{\partial \theta _{i}}\frac{\partial B^{T}}{\partial \theta j}%
(z^{-1}I-A)^{-T}C^{T}\right] z^{-1}dz$

\bigskip

Case 2: \ \ i,j$\in $I

\bigskip

$g_{ij}=\frac{1}{2\pi \text{i}}\doint_{C}tr$

$\left[ C(zI-A)^{-1}\frac{\partial A}{\partial \theta _{i}}%
(zI-A)^{-1}(z^{-1}I-A)^{-T}\frac{\partial A^{T}}{\partial \theta _{j}}%
(z^{-1}I-A)^{-T}C^{T}\right] z^{-1}dz$

\bigskip

Case 3: \ \ \ i$\in $I \ \ and \ \ j$\in $J \ \ 

\bigskip

$g_{ij}=$

$=\frac{1}{2\pi \text{i}}\doint_{C}tr\left[ C(zI-A)^{-1}\frac{\partial A}{%
\partial \theta _{i}}(zI-A)^{-1}\frac{\partial B^{T}}{\partial \theta j}%
(z^{-1}I-A)^{-T}C^{T}\right] z^{-1}dz$

\begin{remark}
\bigskip In the deterministic case, letting $\theta _{i}=\left[ K\right]
_{kl},i=nm+1,...,2nm,$ with k = integer(i/m)-n+1 and l = i mod(m), which
amounts to naming $\theta _{i}$ in lexicographic order in matrix K, one has $%
\frac{\partial K}{\partial \theta _{i}}\frac{\partial K^{T}}{\partial \theta
_{j}}=\mathbb{O}$ (null matrix) whenever i mod(m) $\neq $ j mod(m) (x mod(y)
meaning the rest of division of x by y). Thus, the display of the metric
tensor formulas for the ARMA and state space representations shows
that g$_{ij}=0$ for certain pairs (i,j), since in case 1 the expression $%
\frac{\partial K}{\partial \theta _{i}}\frac{\partial K^{T}}{\partial \theta
_{j}}$(=$\frac{\partial B}{\partial \theta _{i}}\frac{\partial B^{T}}{%
\partial \theta j}$, as B=K) appears in both representations.
\end{remark}

\subsection{Stochastic Metric Tensor}

Stochastic linear dynamical systems are here defined as those ones which,
besides having an input channel for known vector sequences, are permanently
being excited by a vector white noise, of which only the first and second
moments are known, the sequence itself being unknown. Hanzon proposes an
internal product for the tangent space of the manifold of these systems
without known input, of fixed McMillan degree and equal number of components
for the vectors of white noise input and colored noise output. The
corresponding tensor formula is obtained here, based upon theorem 5.1.1. A
simpler formula is also presented, derived from another suggested internal
product by that author. Detailed calculations for their applications to some
simple examples are shown.

\bigskip

A stationary time series will be considered here as a realization of a
gaussian and stationary stochastic process \{y$_{t}\}$ . It will be supposed
that y has zero mean, y $\in \mathcal{R}^{m}$ and t =1,2,3,... , that is,
time will be discrete. Thus, \{y$_{t}\}$ is totally characterized by its
auto-covariance function

\begin{equation}
\Gamma _{i}=E(y_{t}y_{t+i}^{T}),i=0,1,2,...,
\end{equation}

where E(.) is the expectance operator over y 's probability distribution
function and the superscript T stands for transposition. This stems from the
fact that the joint probability density of a string of size I+1 of the
process, that is, of the random vector w$%
^{T}=(y_{t}^{T},y_{t+1}^{T},...,y_{t+I}^{T})$ is given by (\cite{GP}, pg.
90):

\begin{equation}
f(w)=[(2\pi )^{m(I+1)}\det G]^{-1/2}e^{-\frac{1}{2}w^{T}Gw}
\end{equation}

where

\begin{equation}
G=\left[ 
\begin{array}{lllll}
\Gamma _{0} & \Gamma _{1} & \Gamma _{2} & ... & \Gamma _{I} \\ 
\Gamma _{1}^{T} & \Gamma _{0} & \Gamma _{1} & ... & \Gamma _{I-1} \\ 
... & ... & ... & ... & ... \\ 
\Gamma _{I}^{T} & \Gamma _{I-1}^{T} & \Gamma _{I-2}^{T} & ... & \Gamma _{0}
\end{array}
\right] \newline
\end{equation}

Being \{y\} stationary, it follows that $\Gamma _{i}=\Gamma _{-i}^{T}.$

The density above presupposes that G be nonsingular, which is equivalent to
w's components being linearly independent.

There is a biunivocal relation between the autocovariance function of a
stochastic process and its spectral density \cite{DH}. In this case, calling
T(z) the spectral density, one has (\cite{AK}, pg. 69):

\begin{equation}
T(z)=\sum_{i=-\infty }^{\infty }\Gamma _{i}z^{-i}
\end{equation}

Now, T(z) can always be factored as T(z)=H(z)RH$^{T}$(z$^{-1}$), where H(z)
is a rational and stable matrix, with stable inverse, and R is symmetric
positive definite.

Imposing, further, that H(z) be causal and H$_{0}$=I, there is a unique pair
[H(z),R] correspondent to rational T(z) (\cite{GP}, pg. 72).

H(z) may be interpreted as the transfer function of a linear system. As a
result, \{y\} will be interpreted as the output of a stable linear system
whose input is a unobserved white noise.

\bigskip

\bigskip

Hanzon \cite{HZ} proposes an internal product for the manifold of stochastic
systems with a common McMillan degree analogous to the deterministic case:

\begin{equation}
<s_{1},s_{2}>=tr\left[ \dsum\limits_{i\text{=}0}^{\infty }\Gamma
_{i}^{1}(\Gamma _{i}^{2})^{T}\right]
\end{equation}

where s$_{1}$ e s$_{2}$ are systems belonging to the tangent bundle of the
manifold.

This is a metric of the covariance system, whose Markov parameters are $%
\{\Gamma _{i}\}.$ So the immediate extension of theorem 5.1.1 to the
stochastic case is:

\begin{theorem}
\textbf{\ }A stochastic Riemannian metric tensor can be obtained by the
formula
\end{theorem}

\begin{equation}
g_{ij}=\frac{1}{2\pi \text{i}}\doint_{C}tr\left[ \frac{\partial U(z)}{%
\partial \theta _{i}}\frac{\partial U^{T}(z^{-1})}{\partial \theta _{j}}%
\right] z^{-1}dz  \label{gM}
\end{equation}

where:

\begin{equation}
U(z)=\sum_{i=0}^{\infty }\Gamma _{i}z^{-i},
\end{equation}

$\Gamma _{i}$= E(y$_{t}$y$_{t+i}^{T}),i=0,1,2,...,$ are the covariances of
the stochastic process\ \{y$_{t}$\}generated by the stochastic system [$%
\{H_{i}\},R$] ,

tr stands for trace,

the superscript T indicates matrix transposition,

''i '' is $\sqrt{-1\text{ }}$and

C is the unitary circle centered in the origin of the complex plane.

A more convenient metric, which can be expressed directly in terms of H(z)
and R, is the one induced by the internal product defined over the two-sided
infinite sequence $\left\{ \Gamma _{1},i=...,-2,-1,0,1,2,...\right\} $
(reminding that $\Gamma _{i}=\Gamma _{-i}^{T},\forall i\in \Bbb{Z}$):

\begin{equation}
<s_{1},s_{2}>=tr\left[ \dsum\limits_{i\text{=}-\infty }^{\infty }\Gamma
_{i}^{1}(\Gamma _{i}^{2})^{T}\right]
\end{equation}

Hanzon (\cite{HZ}, pg. 208) refers to this choice as \textit{an also quite
attractive possibility. }Defining T(z) = $\sum_{i=-\infty }^{\infty }\Gamma
_{i}z^{-i}$, it follows the \bigskip

\begin{theorem}
\textbf{\ }A stochastic Riemannian metric tensor can be obtained by the
formula
\end{theorem}

\begin{equation}
g_{ij}=\frac{1}{2\pi \text{i}}\doint_{C}tr\left[ \frac{\partial T(z)}{%
\partial \theta _{i}}\frac{\partial T^{T}(z^{-1})}{\partial \theta _{j}}%
\right] z^{-1}dz  \label{gT}
\end{equation}

where T(z) = $\sum_{i=-\infty }^{\infty }\Gamma _{i}z^{-i}$ = H(z)RH$^{T}$(z$%
^{-1}$) is the spectral density of the output $\left\{ y_{t}\right\} $ of
the linear dynamic system whose transfer function is H(z) and whose input is
a white noise $\left\{ \varepsilon _{t}\right\} $ of covariance matrix R and
zero mean.

\bigskip

\textbf{Proof: } Percival's formula

\begin{equation}
tr\left[ \dsum\limits_{i=0}^{\infty }A_{i}^{1}(A_{i}^{2})^{T}\right] =\frac{1%
}{2\pi \text{i}}\doint_{C}tr\left[ A_{1}(z)A_{2}^{T}(z^{-1})\right] z^{-1}dz
\end{equation}

which was the kernel of the proof of theorem 5.1.1, relates a summation of
the sequence $\left\{ A_{i},i\in \Bbb{N}\right\} $ to an integral of the
z-transform A(z) = $\sum\limits_{i=0}^{\infty }A_{i}z^{-1}$. It can also be
stated as

\begin{equation}
tr\left[ \dsum\limits_{i=-\infty }^{\infty }A_{i}^{1}(A_{i}^{2})^{T}\right] =%
\frac{1}{2\pi \text{i}}\doint_{C}tr\left[ A_{1}(z)A_{2}^{T}(z^{-1})\right]
z^{-1}dz
\end{equation}

where A(z) is now defined by A(z) = $\sum\limits_{i=-\infty }^{\infty
}A_{i}z^{-1}$, relating, thus, the summation of the sequence $\left\{
A_{i},i\in \Bbb{Z}\right\} $ to an integral of this last z-transform.

The remainder of the proof is equal to that of theorem 5.1.

\textbf{Q.E.D.}

\bigskip

The relation between T(z) and U(z) is the following:

\begin{equation}
T(z)=U(z)+U^{T}(z^{-1})-\Gamma _{0}  \label{TM}
\end{equation}

\textit{\bigskip }

\subsection{Some examples of tensors for the stochastic case}

\bigskip

Consider, as examples, the following stochastic systems\ (all of them
scalar, that is, $y_{t},x_{t},\varepsilon _{t}\in \Bbb{R},$ and, to
simplify, R=1):

\subsubsection{ In the state space representation}

\bigskip

\begin{equation}
x_{t+1}=Ax_{t}+B\varepsilon _{t}
\end{equation}

\begin{equation}
y_{t}=Cx_{t}+D\varepsilon _{t}
\end{equation}

\bigskip

$y_{t},x_{t},\varepsilon _{t}\in \Bbb{R}.$

In this case, the covariance R of $\varepsilon _{t}$ is a scalar and has a
purely multiplicative effect on the tensor.

\bigskip

\textbf{Example 1}

\bigskip

To simplify, let R=1, D=0 and C=1. Then, the system reduces to

\begin{equation}
x_{t+1}=ax_{t}+b\varepsilon _{t}
\end{equation}

\begin{equation}
y_{t}=x_{t}
\end{equation}

whose transfer function is

\begin{equation}
h(z)=\frac{b}{z-a}
\end{equation}

So, its spectral density is

\begin{equation}
T(z)=H(z)RH^{T}(z^{-1})=\frac{b}{z-a}\frac{b}{z^{-1}-a}=\allowbreak -b^{2}%
\frac{z}{\left( -a+z\right) \left( -1+az\right) }  \label{T1}
\end{equation}

The partial derivatives become

\bigskip

\bigskip $\frac{\partial T(z)}{\partial a}=\allowbreak \allowbreak b^{2}z%
\frac{1-2az+z^{2}}{\left( -a+z\right) ^{2}\left( -1+az\right) ^{2}}$

$\frac{\partial T(z)}{\partial b}=\allowbreak -2b\frac{z}{\left( z-a\right)
\left( -1+az\right) }$

\bigskip

Then, noticing that T(z)=T(z$^{-1}$),

$g_{11}=\frac{1}{2\pi \text{i}}\doint_{C}tr\left[ \frac{\partial T(z)}{%
\partial a}\frac{\partial T^{T}(z^{-1})}{\partial a}\right] z^{-1}dz=$

=$\frac{1}{2\pi \text{i}}\doint_{C}\left[ \allowbreak b^{2}z\frac{1-2az+z^{2}%
}{\left( -a+z\right) ^{2}\left( -1+az\right) ^{2}}\right] ^{2}z^{-1}dz$

$g_{22}=\frac{1}{2\pi \text{i}}\doint_{C}tr\left[ \frac{\partial T(z)}{%
\partial b}\frac{\partial T^{T}(z^{-1})}{\partial b}\right] z^{-1}dz=$

=$\frac{1}{2\pi \text{i}}\doint_{C}\left[ -2b\frac{z}{\left( z-a\right)
\left( -1+az\right) }\right] ^{2}z^{-1}dz$

$g_{12}=g_{21}=\frac{1}{2\pi \text{i}}\doint_{C}tr\left[ \frac{\partial T(z)%
}{\partial a}\frac{\partial T^{T}(z^{-1})}{\partial b}\right] z^{-1}dz$

$g_{12}=g_{21}=$

=$\frac{1}{2\pi \text{i}}\doint_{C}\left[ -2b\frac{z}{\left( z-a\right)
\left( -1+az\right) }\right] \left[ b^{2}z\frac{1-2az+z^{2}}{\left(
z-a\right) ^{2}\left( -1+az\right) ^{2}}\right] z^{-1}dz$

\bigskip

Calculating the integrals, the following tensor is obtained:

\begin{equation}
G=\left[ 
\begin{array}{ll}
-2\dfrac{(2a^{4}+7a^{2}+1)b^{4}}{(a^{2}-1)^{5}} & 4\dfrac{ab^{3}(a^{2}+2)}{%
(a^{2}-1)^{4}} \\ 
4\dfrac{ab^{3}(a^{2}+2)}{(a^{2}-1)^{4}} & -4\dfrac{b^{2}(a^{2}+1)}{%
(a^{2}-1)^{3}}
\end{array}
\right]  \label{G1}
\end{equation}

This is an interesting example, because it is also exhibited in \cite{HZ},
pg. 223, using the metric of theorem 5.2 (but through another method,
involving Ricati and Lyapunov equations and not a complex integral. The
tensor obtained in that work is: 
\begin{equation}
G=\left[ 
\begin{array}{ll}
-\dfrac{(9a^{2}+1)b^{4}}{(a^{2}-1)^{5}} & 6\dfrac{ab^{3}}{(a^{2}-1)^{4}} \\ 
6\dfrac{ab^{3}}{(a^{2}-1)^{4}} & -4\dfrac{b^{2}}{(a^{2}-1)^{3}}
\end{array}
\right]  \label{G2}
\end{equation}

Notice that the denominators coincide with those of the tensor($\ref{G1}$) ,
but not the numerators. The tensor (\ref{G2}) can be found by theorem 5.2 by
the following process:

Expand the transfer function in power series

\begin{equation}
h(z)=\frac{b}{z-a}=\allowbreak
bz^{-1}(1+az^{-1}+a^{2}z^{-2}+...)=bz^{-1}+abz^{-2}+a^{2}bz^{-3}+...
\end{equation}

to obtain the system's Markov parameters

\begin{equation}
h_{0}=0,h_{1}=b,h_{2}=ab,h_{3}=a^{2}b,...,h_{i}=a^{i-1}b
\end{equation}

The sequence of covariances $\left\{ \Gamma _{k}\right\} $ can be found by:

\bigskip

\bigskip $\Gamma _{k}=E(y_{t}y_{t+k}^{T})=E\left[ \sum\limits_{i=0}^{\infty
}H_{i}\varepsilon _{t-i}\left( \sum\limits_{j=0}^{\infty }H_{j}\varepsilon
_{t+k-j}\right) ^{T}\right] =$

=$E\left[ \sum\limits_{i=0}^{\infty }H_{i}\varepsilon
_{t-i}\sum\limits_{j=0}^{\infty }\varepsilon _{t+k-j}^{T}H_{j}^{T}\right] $

$\Gamma _{k}=\sum\limits_{i=0}^{\infty }\sum\limits_{j=0}^{\infty }H_{i}E%
\left[ \varepsilon _{t-i}\varepsilon _{t+k-j}^{T}\right] H_{j}^{T}=\sum%
\limits_{i=0}^{\infty }\sum\limits_{j=0}^{\infty }H_{i}R\delta
_{j-k,i}H_{j}^{T}=$

$=\sum\limits_{i=0}^{\infty }H_{i}RH_{i+k}^{T}$

where $\delta _{j-k,i}$ is Kronecker's delta and reminding that y$%
_{t}=\sum\limits_{i=0}^{\infty }H_{i}\varepsilon _{t-i}.$

Substituting the Markov sequence of the example, comes

\begin{equation}
\Gamma _{k}=\sum\limits_{i=1}^{\infty }a^{i-1}b^{2}a^{i+k-1}
\end{equation}

Then, 
\begin{equation}
U(z)=\sum\limits_{k=0}^{\infty }(\sum\limits_{i=1}^{\infty
}a^{i-1}b^{2}a^{i+k-1})z^{-k}=\allowbreak -b^{2}\frac{z}{\left(
a^{2}-1\right) \left( -a+z\right) }
\end{equation}

Substituting in formula (\ref{gM}) of theorem 5.2, tensor (\ref{G2}) is
obtained.

\medskip

Notice that equation (\ref{TM}) is satisfied:

\bigskip

$U(z^{-1})=-\frac{b^{2}}{z\left( a^{2}-1\right) \left( -a+\frac{1}{z}\right) 
}$

\begin{equation}
\Gamma _{0}=\sum\limits_{i=0}^{\infty
}H_{i}RH_{i}^{T}=\sum\limits_{i=1}^{\infty }a^{i-1}b^{2}a^{i-1}=-\frac{b^{2}%
}{a^{2}-1}
\end{equation}

$U(z)+U^{T}(z^{-1})-\Gamma _{0}=$

$=-b^{2}\frac{z}{\left( a^{2}-1\right) \left( -a+z\right) }+\left( -\frac{%
b^{2}}{z\left( a^{2}-1\right) \left( -a+\frac{1}{z}\right) }\right) -\left( -%
\frac{b^{2}}{a^{2}-1}\right) =\allowbreak $

$=-z\frac{b^{2}}{\left( az-1\right) \left( -a+z\right) }=T(z)$ (see
expression (\ref{T1}).

\bigskip

\textbf{Example 2}

\bigskip

\smallskip

The next example is of the important innovations model:

\begin{equation}
x_{t+1}=ax_{t}+b\varepsilon _{t}
\end{equation}

\begin{equation}
y_{t}=x_{t}+\varepsilon _{t}
\end{equation}

The transfer function in this case is

\begin{equation}
h(z)=\frac{z+b-a}{z-a}
\end{equation}

The spectral density becomes

\begin{equation}
T(z)=h(z)h(z^{-1})=\allowbreak \allowbreak \left( z+b-a\right) \frac{-1-bz+az%
}{\left( z-a\right) \left( -1+az\right) }
\end{equation}

The result found by theorem 3 is:

\[
g_{11}=-2\frac{%
(-a^{4}+7b^{2}a^{2}+1+a^{6}-a^{2}-4ba^{5}+2b^{2}a^{4}+4ba+b^{2})b^{2}}{%
(a^{2}-1)^{5}} 
\]

\begin{equation}
g_{12}=g_{21}=2\frac{b(a-2a^{3}+a^{5}-4ba^{4}+2b^{2}a^{3}+3ba^{2}+4b^{2}a+b)%
}{(a^{2}-1)^{4}}
\end{equation}

$g_{22}=-2\dfrac{a^{4}-2a^{2}+1-4ba^{3}+2b^{2}a^{2}+4ba+2b^{2}}{(a^{2}-1)^{3}%
}$

\bigskip

To use theorem 5.2 semi-infinite metric, compute:

$h(z)=\frac{z+b-a}{z-a}=\frac{z-a}{z-a}+\frac{b}{z-a}=1+\frac{b}{z-a}%
=1+bz^{-1}+abz^{-2}+a^{2}bz^{-3}+...$

Thus,

\begin{equation}
h_{i}=a^{i-1}b,i=1,2,3,...
\end{equation}

\begin{equation}
h_{0}=1
\end{equation}

The general expressions

\bigskip

$\Gamma _{k}=\sum\limits_{i=0}^{\infty }H_{i}RH_{i+k}^{T}$

$U(z)=\sum_{k=0}^{\infty }\Gamma _{k}z^{-k}$

\bigskip

become, for this case:

\bigskip

For $k=1,2,3,..$

$\Gamma _{k}=\sum\limits_{i=0}^{\infty
}H_{i}RH_{i+k}^{T}=\sum\limits_{i=1}^{\infty
}H_{i}RH_{i+k}^{T}+H_{0}RH_{k}^{T}$

\begin{equation}
\Gamma _{k}=\sum\limits_{i=1}^{\infty
}a^{i-1}(b)^{2}a^{i+k-1}+a^{k-1}b=\allowbreak -b\frac{a^{k}b-a^{1+k}+a^{k-1}%
}{a^{2}-1}
\end{equation}

\begin{equation}
\Gamma _{0}=\sum\limits_{i=1}^{\infty
}H_{i}RH_{i}^{T}+H_{0}RH_{0}^{T}=\left( \sum\limits_{i=1}^{\infty
}a^{i-1}(b)^{2}a^{i-1}\right) +1=\allowbreak \frac{-b^{2}+a^{2}-1}{a^{2}-1}
\end{equation}

Then, $U(z)=\sum_{k=0}^{\infty }\Gamma _{k}z^{-k}=\sum_{k=1}^{\infty }\Gamma
_{k}z^{-k}+\Gamma _{0}$

$U(z)=\sum\limits_{k=1}^{\infty }-b\frac{a^{k}b-a^{1+k}+a^{k-1}}{a^{2}-1}%
z^{-k}+\allowbreak \frac{-b^{2}+a^{2}-1}{a^{2}-1}=\allowbreak \frac{%
ba^{2}-b-b^{2}z-a^{3}+a^{2}z+a-z}{\left( a^{2}-1\right) \left( -a+z\right) }$

$\allowbreak $

Having U(z), it is enough to apply theorem 5.2 formula (\ref{gM}), to get to
the corresponding tensor:

\[
g_{11}=\allowbreak -\left(
a^{6}-4a^{5}b-a^{4}+9b^{2}a^{2}-a^{2}+4ba+1+b^{2}\right) \frac{b^{2}}{\left(
a-1\right) ^{5}\left( a+1\right) ^{5}} 
\]

\begin{equation}
g_{12}=\allowbreak \left( a^{5}-4ba^{4}-2a^{3}+3ba^{2}+a+6ab^{2}+b\right) 
\frac{b}{\left( a-1\right) ^{4}\left( a+1\right) ^{4}}
\end{equation}

$g_{22}=\allowbreak -\dfrac{-2a^{2}+4b^{2}+1+a^{4}+4ba-4a^{3}b}{\left(
a-1\right) ^{3}\left( a+1\right) ^{3}}$

\bigskip

\subsubsection{ In the ARMA representation}

\bigskip The ARMA equations is:

\begin{equation}
A_{0}y_{t}+A_{1}y_{t-1}+...+A_{p}y_{t-p}=B_{0}u_{t}+B_{1}u_{t-1}+...+B_{p}u_{t-p}
\end{equation}

The univariate ARMA(1,1) case will be exhibited here. It is defined by

\bigskip

\bigskip \textbf{Example 3}

\begin{equation}
y_{t}+ay_{t-1}=\varepsilon _{t}+b\varepsilon _{t-1}
\end{equation}

whose transfer function is

\begin{equation}
h(z)=\frac{z+b}{z+a}
\end{equation}

The spectral density becomes

\begin{equation}
T(z)=h(z)h(z^{-1})=\allowbreak \frac{z+b}{z+a}\frac{\frac{1}{z}+b}{\frac{1}{z%
}+a}=\allowbreak \left( z+b\right) \frac{1+bz}{\left( az+1\right) \left(
a+z\right) }
\end{equation}

The following tensor is obtained by theorem 3:

\bigskip

$g_{11}=2\frac{%
-43b^{2}a^{2}-5b^{2}+b^{2}a^{6}-13b^{2}a^{4}-2a^{4}-2b^{4}a^{4}+24b^{3}a^{3}-7a^{2}-7b^{4}a^{2}+16ba+16b^{3}a-1-b^{4}%
}{\left( a^{2}-1\right) ^{5}}$

\begin{equation}
g_{12}=-2\frac{%
-11ba-8ba^{3}+ba^{5}-2b^{3}a^{3}+15b^{2}a^{2}+5a^{2}-a^{2}b^{3}+1+3b^{2}}{%
\left( a^{2}-1\right) ^{4}}
\end{equation}

$g_{22}=2\dfrac{a^{4}-4a^{2}-1-2b^{2}a^{2}+8ba-2b^{2}}{(a^{2}-1)^{3}}$

\bigskip

The computations to obtain the tensor prescribed by theorem 5.2 are:

\bigskip

$h(z)=\frac{z+b}{z+a}=\frac{z}{z+a}+\frac{b}{z+a}$

\bigskip

$\frac{z}{z+a}=\allowbreak \frac{1}{1+az^{-1}}%
=1-az^{-1}+a^{2}z^{-2}-a^{3}z^{-3}+...$

\bigskip

$\frac{b}{z+a}=bz^{-1}\allowbreak \frac{1}{1+az^{-1}}%
=bz^{-1}-baz^{-2}+ba^{3}z^{-3}+...$

\bigskip

Thus,

\begin{equation}
h(z)=1-(a-b)z^{-1}+a(a-b)z^{-2}-a^{2}(a-b)z^{-3}+...
\end{equation}

So,

\begin{equation}
h_{i}=a^{i-1}(a-b)(-1)^{i},i=1,2,3,...
\end{equation}

\begin{equation}
h_{0}=1
\end{equation}

$\allowbreak $Now, for $k=1,2,3,..$

\bigskip

$\Gamma _{k}=\sum\limits_{i=1}^{\infty
}a^{i-1}(a-b)(-1)^{i}a^{i+k-1}(a-b)(-1)^{i+k}+a^{k-1}(a-b)(-1)^{k}%
\allowbreak $%
\begin{equation}
\Gamma _{k}=\left( a-b\right) \left( -1\right) ^{k}\frac{a^{k}b-a^{k-1}}{%
a^{2}-1}
\end{equation}

\begin{equation}
\Gamma _{0}=\left( \sum\limits_{i=1}^{\infty
}a^{i-1}(a-b)^{2}a^{i-1}(-1)^{2i-2}\right) +1=\allowbreak \frac{2ba-b^{2}-1}{%
a^{2}-1}
\end{equation}

Then,

$U(z)=\sum\limits_{k=1}^{\infty }\left( a-b\right) \left( -1\right) ^{k}%
\frac{a^{k}b-a^{k-1}}{a^{2}-1}z^{-k}+\allowbreak \allowbreak \frac{%
2ba-b^{2}-1}{a^{2}-1}=\allowbreak \frac{ba^{2}-b+2baz-b^{2}z-z}{\left(
a^{2}-1\right) \left( a+z\right) }$

The corresponding tensor is:

\bigskip

$g_{11}=$

$=\allowbreak \allowbreak \frac{%
3b^{2}a^{6}-4b^{3}a^{5}-4a^{5}b-7b^{2}a^{4}+24a^{3}b+24b^{3}a^{3}-49b^{2}a^{2}-9a^{2}-9a^{2}b^{4}+20b^{3}a+20ba-1-b^{4}-7b^{2}%
}{\left( a-1\right) ^{5}\left( a+1\right) ^{5}}\qquad \qquad $

$g_{12}=\allowbreak -\frac{%
3a^{5}b-2a^{4}-4b^{2}a^{4}-6a^{3}b+7a^{2}+17b^{2}a^{2}-15ba-6b^{3}a+1+5b^{2}%
}{\left( a-1\right) ^{4}\left( a+1\right) ^{4}}$

\begin{equation}
g_{22}=\allowbreak \dfrac{3a^{4}-4a^{3}b-6a^{2}+12ba-1-4b^{2}}{\left(
a-1\right) ^{3}\left( a+1\right) ^{3}}
\end{equation}
\qquad \qquad \qquad

\bigskip

\section{Conclusions}

The differentiable manifold of systems of a common McMillan degree is a
natural set for the state space representation, because this degree is the
number of components of the state in the minimal representation, thus
appearing explicitly in the models. It gives rise to the so-called
overlapping parametrizations, which exhibit much more flexibility than the
canonical ones, which is important for the numerical process of
identification, given the possibility of changing of model even on-line,
whenever a malconditioning is detected.

For the ARMA representation, the establishing of overlapping

parametrizations based upon classes of models representing systems with a
common McMillan degree results in a clumsy class of mathematical models \cite
{CG}, because this degree is not natural for that representation.
Unfortunately, the natural integer - p, the minimum degree of the AR and MA
polynomials - doesn't give rise, as shown, to a class of models whose image
(the systems of a common p) is a differentiable manifold, thus not allowing
the definition of a natural ARMA overlapping parametrization.

A remedial solution could be to work with a rougher class of models (the
ARMA(p,p)-irreducible ones), which would not be strictly identifiable for
systems in which not all Kronecker indices are equal; since the set of
systems in which all them are equal is generic, maybe this was not such a
big handicap...

The complexity of the Riemannian metric tensor formulas grow exponentially
with the increase of systems dimensions. There are two ways of attenuating
this problem:

1) Computer languages for algebraic symbolic processing like MAPLE and
MATHEMATICA$^{TM}$ are sufficiently flexible to integrate numerical
algorithms with symbolic ones, thus dispensing the need of manual
transcription of formulas.

2) Since the only systems treated in the theory here presented are the
stable ones, the parameters that appear in the denominator of the transfer
function are bounded by restriction such as being between -1 and 1. In this
case, certain terms of higher order in the expressions obtained for the
metric tensors can be despised without considerable loss of precision, thus
reducing the complexity of the formulas.

The availability of powerful computer packages for algebraic computation \
turns attractive the use of analytical formulas involving complex integrals
and partial derivatives of polynomial matrices for the study of geometrical
properties of spaces of linear dynamical systems.

Given the close relations between the Riemannian metric tensor, the Fisher
information matrix, the covariance matrix of the parameters estimators and
the Hessian matrix of some common objective functions used in parametric
identification, the author believes that the results displayed in this
article hold some relevance for the classical problem of linear dynamical
systems identification.


\bigskip

\begin{thebibliography}{99}
\bibitem{AK}  Aoki, M. (1987). \textbf{State space modeling of time series}.
Springer Verlag.

\bibitem{AP}  Aplevich, J.D. \textbf{Singular Pencil Models in Systems
Design and Control}, Internal Report N2L 3G1, Elet. Eng. Dpt., University
of. Waterloo, 1981.

\bibitem{BJ}  Box, G.E.P. and Jenkins,G.M. \textbf{Time Series Analysis,
Forecasting and Control}. Holden-Day. San Francisco. 1970.

\bibitem{CG}  Correa, G. O. and Glover, K. \textbf{Pseudo-canonical forms,
identifiable parametrizations and simple estimation for linear multivariable
systems: input-output models.} Automatica, vol. 20, n.4, pp. 429-442. 1984.

\bibitem{CH}  Chou, T. C., \textbf{Geometry of Linear Systems and
Identification}. Ph.D. thesis. Trinity College, Cambridge,1994.

\bibitem{CK}  Clark, J.M.C.\textbf{The consistent selection of local
coordinates in linear systems identification},\ JAAC Purdue University,
Lafayette, Indiana, 1976, 576-580.

\bibitem{DH}  Denham, M. J. (1974). \textbf{Canonical forms for the
identification of multivariable linear systems}. IEEE Trans. Autom. Control.
Vol AC-19, n.6, pp. 646-656.

\bibitem{GD}  Guidorzi, R.P. \textbf{Invariants and canonical forms for
systems structural and parametric identification}. Automatica, vol.17,
n.1,pp117-133.

\bibitem{GP}  Goodwin, G. C. e Payne, R.L. (1977).\textbf{\ Dynamic system
identification}. Academic Press.

\bibitem{GW}  Gevers, M. and Wertz,V. \textbf{\ Uniquely identifiable
state-space and ARMA parametrizations for multivariable systems},
Automatica, vol. 20, n. 3, 1984, 333-347.

\bibitem{HD}  Hannan, E.J. and Deistler, M. \textbf{The statistical theory
of linear systems}. John Wiley and Sons. 1988.

\bibitem{HZ}  Hanzon, B. \textbf{Identifiability, recursive identification
and spaces of linear dynamical systems}, Ph.D. Thesis, Department of
Econometrics, Erasmus University, Rotterdam, 1986.

\bibitem{KT}  Kailath, T. \textbf{Linear Systems}. Prentice Hall.1980.

\bibitem{OP}  Oppenheim, A.V. and Shafer,\ R.W. \textbf{Digital signal
processing}, Prentice-Hall, 1976, 66.

\bibitem{PT}  Peeters, R. \textbf{System identification based on Riemannian
geometry: theory and algorithms}, Ph.D. Thesis, Free University of
Amsterdam, 1994, and Research Report nr. 64, Tinbergen Institute Research
Series, Tinbergen Institute, Rotterdam.

\bibitem{TT}  Tiao, G.C. and Tsay, R.S. \textbf{Multiple time series
modeling and extended sample cross correlations}. Technical report n. 690.
Statistics Dept. Univesity of Wisconsin. 1982.
\end{thebibliography}
\end{document}